\documentclass[b5paper,reqno]{amsart}%
\usepackage{amsfonts}
\usepackage{amsmath}
\usepackage{amssymb}
\usepackage{amsthm}
\usepackage{enumerate}
\usepackage[mathscr]{eucal}
\usepackage{graphicx}%
\newtheorem{theorem}{Theorem}
\theoremstyle{plain}

\newtheorem{corollary}{Corollary}

\newtheorem{definition}{Definition}
\newtheorem{example}{Example}

\newtheorem{lemma}{Lemma}

\newtheorem{remark}{Remark}

\numberwithin{equation}{section}
\begin{document}
\title[Periodicity on time scales]{A new periodicity concept for time scales}
\author{Murat Ad\i var}
\newcommand{\acr}{\newline\indent}
\address{\llap{*\,}Izmir University of Economics \acr
                Department of Mathematics \acr
                35330, Bal\c{c}ova\acr
                Izmir, Turkey}
\email{murat.adivar@ieu.edu.tr}
\urladdr{http://homes.ieu.edu.tr/\symbol{126}madivar}
\date{April 19, 2011}

\begin{abstract}
By means of the shift operators we introduce a new periodicity concept on time
scales. This new approach will enable researchers to investigate periodicity
notion on a large class of time scales whose members may not satisfy the
condition%
\[
\text{"there exists a }P>0\text{ such that }t\pm P\in\mathbb{T}\text{ for all
}t\in\mathbb{T}\text{",}%
\]
which is being currently used. Therefore, the results of this paper open an
avenue for the investigation of periodic solutions of $q$-difference equations
and more.

\end{abstract}
\maketitle

%\address{Izmir University of Economics, Department of Mathematics\\
%35330, Bal\c{c}ova, \.{I}zmir, Turkey}

\section{Introduction}

In the last two decades, theory of time scales has become a very useful tool
for the unification of difference and differential equations under dynamic
equations on time scales (see \cite{adivar}-\cite{Kaufmann&raffoul}, and
references therein). A time scale, denoted by $\mathbb{T}$, is a non-empty
arbitrary subset of real numbers. To be able to investigate the notion of
periodicity of the solutions of dynamic equations on time scales researchers
had to first introduce the concept of periodic time scales and then define
what it meant for a function to be periodic on such a time scale. To be more
specific, we restate the following definitions and introductory examples which
can be found in \cite{atici}, \cite{Bi&Bohner}, and \cite{Kaufmann&raffoul}.

\begin{definition}
[9]\label{oldperts}A time scale $\mathbb{T}$ is said to be \emph{periodic} if
there exists a $P>0$ such that $t\pm P\in\mathbb{T}$ for all $t\in\mathbb{T}$.
If $\mathbb{T}\neq\mathbb{R}$, the smallest positive $P$ is called the
\emph{period} of the time scale.
\end{definition}

\begin{example}
The following time scales are periodic.

\begin{enumerate}
\item[i.] $\mathbb{T}=\mathbb{Z}$ has period $P=1$,

\item[ii.] $\mathbb{T}=h\mathbb{Z}$ has period $P=h$,

\item[iii.] $\mathbb{T}=\mathbb{R}$,

\item[iv.] $\mathbb{T}=\bigcup_{i=-\infty}^{\infty}[(2i-1)h,2ih],\,h>0$ has
period $P=2h$,

\item[v.] $\mathbb{T}=\{t=k-q^{m}:k\in\mathbb{Z},m\in\mathbb{N}_{0}\}$ where,
$0<q<1$ has period $P=1$.
\end{enumerate}
\end{example}

\begin{definition}
\label{old per}Let $\mathbb{T}\neq\mathbb{R}$ be a periodic time scale with
period $P$. We say that the function $f\colon\mathbb{T}\rightarrow\mathbb{R}$
is periodic with period $T$ if there exists a natural number $n$ such that
$T=nP$, $f(t\pm T)=f(t)$ for all $t\in\mathbb{T}$ and $T$ is the smallest
number such that $f(t\pm T)=f(t)$. If $\mathbb{T}=\mathbb{R}$, we say that $f$
is periodic with period $T>0$ if $\,T$ is the smallest positive number such
that $f(t\pm T)=f(t)$ for all $t\in\mathbb{T}$.
\end{definition}

Based on the Definitions \ref{oldperts} and \ref{old per}, periodicity and
existence of periodic solutions of dynamic equations on time scales were
studied by various researchers and for first papers on the subject we refer to
(see for instance \cite{adivar1}-\cite{Bi&Bohner}, \cite{Kaufmann 2}-\cite{li}).

There is no doubt that a time scale $\mathbb{T}$ which is periodic in the
sense of Definition \ref{oldperts} must satisfy%
\begin{equation}
t\pm P\in\mathbb{T}\text{ for all }t\in\mathbb{T} \label{addition}%
\end{equation}
for a fixed $P>0$. This property obliges the time scale to be unbounded from
above and below. However, these two restrictions prevents us from
investigating the periodic solutions of $q$-difference equations since the
time scale
\[
\overline{q^{\mathbb{Z}}}=\left\{  q^{n}:q>1\text{ is constant and }%
n\in\mathbb{Z}\right\}  \cup\left\{  0\right\}
\]
which is neither\ closed under the operation $t\pm P$ for a fixed $P>0$ nor
unbounded below.

The main purpose of this paper is to introduce a new periodicity concept on
time scales which does not oblige the time scale to be closed under the
operation $t\pm P$ for a fixed $P>0$ or to be unbounded. We define our new
periodicity concept with the aid of shift operators which are first defined in
\cite{adivar} and then generalized in \cite{bams}.

\section{Shift operators}

Next, we give a generalized version of shift operators (see \cite{bams}). A
limited version of shift operators can be found in \cite{adivar}. Hereafter,
we use the notation $\left[  a,b\right]  _{\mathbb{T}}$ to indicate the time
scale interval $[a,b]\cap\mathbb{T}$ . The intervals $[a,b)_{\mathbb{T}}$,
$(a,b]_{\mathbb{T}}$, and $\left(  a,b\right)  _{\mathbb{T}}$ are similarly defined.

\begin{definition}
\label{shift} Let $\mathbb{T}^{\ast}$ be a non-empty subset of the time scale
$\mathbb{T}$ including a fixed number $t_{0}\in\mathbb{T}^{\ast}$ such that
there exist operators $\delta_{\pm}:[t_{0},\infty)_{\mathbb{T}}\times
\mathbb{T}^{\ast}\rightarrow\mathbb{T}^{\ast}$ satisfying the following properties:

\begin{enumerate}
\item[P.1] The functions $\delta_{\pm}$ are strictly increasing with respect
to their second arguments, i.e., if
\[
(T_{0},t),(T_{0},u)\in\mathcal{D}_{\pm}:=\left\{  (s,t)\in\lbrack t_{0}%
,\infty)_{\mathbb{T}}\times\mathbb{T}^{\ast}:\delta_{\pm}(s,t)\in
\mathbb{T}^{\ast}\right\}  ,
\]
then
\[
T_{0}\leq t<u\text{ implies }\delta_{\pm}(T_{0},t)<\delta_{\pm}(T_{0},u),
\]

\item[P.2] If $(T_{1},u),(T_{2},u)\in\mathcal{D}_{-}$ with $T_{1}<T_{2}$, then%
\[
\delta_{-}(T_{1},u)>\delta_{-}(T_{2},u),
\]
and if $(T_{1},u),(T_{2},u)\in\mathcal{D}_{+}$ with $T_{1}<T_{2}$, then
\[
\delta_{+}(T_{1},u)<\delta_{+}(T_{2},u),
\]

\item[P.3] If $t\in\lbrack t_{0},\infty)_{\mathbb{T}}$, then $(t,t_{0}%
)\in\mathcal{D}_{+}$ and $\delta_{+}(t,t_{0})=t$. Moreover, if $t\in
\mathbb{T}^{\ast}$, then $(t_{0},t)$ $\in\mathcal{D}_{+}$ and $\delta
_{+}(t_{0},t)=t$ holds,

\item[P.4] If $(s,t)\in\mathcal{D}_{\pm}$, then $(s,\delta_{\pm}%
(s,t))\in\mathcal{D}_{\mp}$ and $\delta_{\mp}(s,\delta_{\pm}(s,t))=t$, respectively,

\item[P.5] If $(s,t)\in\mathcal{D}_{\pm}$ and $(u,\delta_{\pm}(s,t))\in
\mathcal{D}_{\mp}$, then $(s,\delta_{\mp}(u,t))\in\mathcal{D}_{\pm}$ and
$\delta_{\mp}(u,\delta_{\pm}(s,t))=\delta_{\pm}(s,\delta_{\mp}(u,t))$, respectively.
\end{enumerate}

\noindent Then the operators $\delta_{-}$ and $\delta_{+}$ associated with
$t_{0}\in\mathbb{T}^{\ast}$ (called the initial point) are said to be
\textit{backward and forward shift operators} on the set $\mathbb{T}^{\ast}$,
respectively. The variable $s\in\lbrack t_{0},\infty)_{\mathbb{T}}$ in
$\delta_{\pm}(s,t)$ is called the shift size. The values $\delta_{+}(s,t)$ and
$\delta_{-}(s,t)$ in $\mathbb{T}^{\ast}$ indicate $s$ units translation of the
term $t\in\mathbb{T}^{\ast}$ to the right and left, respectively. The sets
$\mathcal{D}_{\pm}$ are the domains of the shift operators $\delta_{\pm}$, respectively.
\end{definition}

Hereafter, we shall denote by $\mathbb{T}^{\ast}$ the largest subset of the
time scale $\mathbb{T}$ such that the shift operators $\delta_{\pm}%
:[t_{0},\infty)_{\mathbb{T}}\times\mathbb{T}^{\ast}\rightarrow\mathbb{T}%
^{\ast}$ exist.

\begin{example}
Let $\mathbb{T=R}$ and $t_{0}=1$. The operators%
\begin{equation}
\delta_{-}(s,t)=\left\{
\begin{array}
[c]{cc}%
t/s & \text{if }t\geq0\\
st & \text{if }t<0
\end{array}
\right.  ,\ \ \ \text{for }s\in\lbrack1,\infty) \label{rs1}%
\end{equation}
and%
\begin{equation}
\delta_{+}(s,t)=\left\{
\begin{array}
[c]{cc}%
st & \text{if }t\geq0\\
t/s & \text{if }t<0
\end{array}
\right.  ,\ \ \ \text{for }s\in\lbrack1,\infty) \label{rs2}%
\end{equation}
are backward and forward shift operators (on the set $\mathbb{R}^{\ast
}=\mathbb{R-}\left\{  0\right\}  $) associated with the initial point
$t_{0}=1$. In the table below, we state different time scales with their
corresponding shift operators.
\[%
\begin{tabular}
[c]{|c||c|c|c|c|}\hline
$\mathbb{T}$ & $t_{0}$ & $\mathbb{T}^{\ast}$ & $\delta_{-}(s,t)$ & $\delta
_{+}(s,t)$\\\hline\hline
$\mathbb{R}$ & $0$ & $\mathbb{R}$ & $t-s$ & $t+s$\\\hline
$\mathbb{Z}$ & $0$ & $\mathbb{Z}$ & $t-s$ & $t+s$\\\hline
$q^{\mathbb{Z}}\cup\left\{  0\right\}  $ & $1$ & $q^{\mathbb{Z}}$ & $\frac
{t}{s}$ & $st$\\\hline
$\mathbb{N}^{1/2}$ & $0$ & $\mathbb{N}^{1/2}$ & $\sqrt{t^{2}-s^{2}}$ &
$\sqrt{t^{2}+s^{2}}$\\\hline
\end{tabular}
\ \ \ \ \ \ \ \ \ \ \
\]

\end{example}

The proof of the next lemma is a direct consequence of Definition \ref{shift}.

\begin{lemma}
\label{lem pro} Let $\delta_{-}$ and $\delta_{+}$ be the shift operators
associated with the initial point $t_{0}$. We have

\begin{enumerate}
\item[i.] $\delta_{-}(t,t)=t_{0}$ for all $t\in\lbrack t_{0},\infty
)_{\mathbb{T}}.$

\item[ii.] $\delta_{-}(t_{0},t)=t$ for all $t\in\mathbb{T}^{\ast},$

\item[iii.] If $(s,t)\in$ $\mathcal{D}_{+}$, then $\delta_{+}(s,t)=u$ implies
$\delta_{-}(s,u)=t$. Conversely, if $(s,u)\in$ $\mathcal{D}_{-}$, then
$\delta_{-}(s,u)=t$ implies $\delta_{+}(s,t)=u$.

\item[iv.] $\delta_{+}(t,\delta_{-}(s,t_{0}))=\delta_{-}(s,t)$ for all
$(s,t)\in$ $\mathcal{D}_{+}$ with $t\geq t_{0,}$

\item[v.] $\delta_{+}(u,t)=\delta_{+}(t,u)$ for all $(u,t)\in\left(  \lbrack
t_{0},\infty)_{\mathbb{T}}\times\lbrack t_{0},\infty)_{\mathbb{T}}\right)
\cap\mathcal{D}_{+}$

\item[vi.] $\delta_{+}(s,t)\in\lbrack t_{0},\infty)_{\mathbb{T}}$ for all
$(s,t)\in$ $\mathcal{D}_{+}$ with $t\geq t_{0,}$,

\item[vii.] $\delta_{-}(s,t)\in\lbrack t_{0},\infty)_{\mathbb{T}}$ for all
$(s,t)\in$ $\left(  [t_{0},\infty)_{\mathbb{T}}\times\lbrack s,\infty
)_{\mathbb{T}}\right)  \cap\mathcal{D}_{-,}$

\item[viii.] If $\delta_{+}(s,.)$ is $\Delta-$differentiable in its second
variable, then $\delta_{+}^{\Delta_{t}}(s,.)>0$,

\item[ix.] $\delta_{+}(\delta_{-}(u,s),\delta_{-}(s,v))=\delta_{-}(u,v)$ for
all $(s,v)\in\left(  \lbrack t_{0},\infty)_{\mathbb{T}}\times\lbrack
s,\infty)_{\mathbb{T}}\right)  \cap\mathcal{D}_{-}$ and $(u,s)\in\left(
\lbrack t_{0},\infty)_{\mathbb{T}}\times\lbrack u,\infty)_{\mathbb{T}}\right)
\cap\mathcal{D}_{-}$,

\item[x.] If $(s,t)\in\mathcal{D}_{-}$ and $\delta_{-}(s,t)=t_{0}$, then $s=t$.
\end{enumerate}
\end{lemma}

\begin{proof}
(i) is obtained from P.3-5 since%
\[
\delta_{-}(t,t)=\delta_{-}(t,\delta_{+}(t,t_{0}))=t_{0}\text{ for all }%
t\in\lbrack t_{0},\infty
)_{\mathbb{T}}.
\]
(ii) is obtained from P.3-P.4 since%
\[
\delta_{-}(t_{0},t)=\delta_{-}(t_{0},\delta_{+}(t_{0},t))=t \text{ for all } t\in\mathbb{T}^{\ast}.
\]
Let $u:=\delta_{+}(s,t)$. By P.4 we have $(s,u)\in\mathcal{D}_{-}$ for all
$(s,t)\in\mathcal{D}_{+}$, and hence,%
\[
\delta_{-}(s,u)=\delta_{-}(s,\delta_{+}(s,t))=t.
\]
The latter part of (iii) can be done similarly. We have (iv) since P.3 and P.5
yield%
\[
\delta_{+}(t,\delta_{-}(s,t_{0}))=\delta_{-}(s,\delta_{+}(t,t_{0}))=\delta
_{-}(s,t).
\]
P.3 and P.5 guarantee that%
\[
t=\delta_{+}(t,t_{0})=\delta_{+}(t,\delta_{-}(u,u))=\delta_{-}(u,\delta
_{+}(t,u))
\]
for all $(u,t)\in\left(  \lbrack t_{0},\infty)_{\mathbb{T}}\times\lbrack
t_{0},\infty)_{\mathbb{T}}\right)  \cap\mathcal{D}_{+}$. Using (iii) we have%
\[
\delta_{+}(u,t)=\delta_{+}(u,\delta_{-}(u,\delta_{+}(t,u)))=\delta_{+}(t,u).
\]
This proves (v). To prove (vi) and (vii) we use P.1-2 to get%
\[
\delta_{+}(s,t)\geq\delta_{+}(t_{0},t)=t\geq t_{0}%
\]
for all $(s,t)\in$ $\left(  [t_{0},\infty)\times\lbrack t_{0},\infty
)_{\mathbb{T}}\right)  \cap\mathcal{D}_{+}$ and%
\[
\delta_{-}(s,t)\geq\delta_{-}(s,s)=t_{0}%
\]
for all $(s,t)\in\left(  \lbrack t_{0},\infty)_{\mathbb{T}}\times\lbrack
s,\infty)_{\mathbb{T}}\right)  \cap\mathcal{D}_{-}$. Since $\delta_{+}(s,t)$
is strictly increasing in its second variable we have (viii) by
\cite[Corollary 1.16]{book2}. (ix) is proven as follows: from P.5 and (v) we
have%
\begin{align*}
\delta_{+}(\delta_{-}(u,s),\delta_{-}(s,v))  &  =\delta_{-}(s,\delta
_{+}(v,\delta_{-}(u,s)))\\
&  =\delta_{-}(s,\delta_{-}(u,\delta_{+}(v,s)))\\
&  =\delta_{-}(s,\delta_{+}(s,\delta_{-}(u,v)))\\
&  =\delta_{-}(u,v)
\end{align*}
for all $(s,v)\in\left(  \lbrack t_{0},\infty)_{\mathbb{T}}\times\lbrack
s,\infty)_{\mathbb{T}}\right)  \cap\mathcal{D}_{-}$ and $(u,s)\in\left(
\lbrack t_{0},\infty)_{\mathbb{T}}\times\lbrack u,\infty)_{\mathbb{T}}\right)
\cap\mathcal{D}_{-}$. Suppose $(s,t)\in\mathcal{D}_{-}$ $=\left\{
(s,t)\in\lbrack t_{0},\infty)_{\mathbb{T}}\times\mathbb{T}^{\ast}:\delta
_{-}(s,t)\in\mathbb{T}^{\ast}\right\}  $ and $\delta_{-}(s,t)=t_{0}$. Then by
P.4 we have%
\[
t=\delta_{+}(s,\delta_{-}(s,t))\in\delta_{+}(s,t_{0})=s.
\]
This is (x). The proof is complete.
\end{proof}

Notice that the shift operators $\delta_{\pm}$ are defined once the initial
point $t_{0}\in\mathbb{T}^{\ast}$ is known. For instance, we choose the
initial point $t_{0}=0$ to define shift operators $\delta_{\pm}(s,t)=t\pm s$
on $\mathbb{T}=\mathbb{R}$. However, if we choose $\lambda\in(0,\infty)$ as
the initial point, then the new shift operators associated with $\lambda$ are
defined by $\widetilde{\delta}_{\pm}(s,t)=t\mp\lambda\pm s$. In terms of
$\delta_{\pm}$ the new shift operators $\widetilde{\delta}_{\pm}$ can be given
as follows%
\[
\widetilde{\delta}_{\pm}(s,t)=\delta_{\mp}(\lambda,\delta_{\pm}(s,t)).
\]

\begin{example}
In the following, we give some particular time scales with shift operators
associated with different initial points to show the change in the formula of
shift operators as the initial point changes.%
\[%
\begin{tabular}
[c]{c||cc|cc|cc}
& \multicolumn{2}{||c|}{$\mathbb{T}=\mathbb{N}^{1/2}$} &
\multicolumn{2}{|c|}{$\mathbb{T}=h\mathbb{Z}$} &
\multicolumn{2}{|c}{$\mathbb{T}=2^{\mathbb{N}}$}\\\hline\hline
$t_{0}$ & $0$ & $\lambda$ & $0$ & $h\lambda$ & $1$ & $2^{\lambda}$\\
$\delta_{-}(s,t)$ & $\sqrt{t^{2}-s^{2}}$ & $\sqrt{t^{2}+\lambda^{2}-s^{2}}$ &
$t-s$ & $t+h\lambda-s$ & $t/s$ & $2^{\lambda}ts^{-1}$\\
$\delta_{+}(s,t)$ & $\sqrt{t^{2}+s^{2}}$ & $\sqrt{t^{2}-\lambda^{2}+s^{2}}$ &
$t+s$ & $t-h\lambda+s$ & $ts$ & $2^{-\lambda}ts$%
\end{tabular}
\ \ \ \ \ \ \
\]
where $\lambda\in\mathbb{Z}_{+}$, $\mathbb{N}^{1/2}=\left\{  \sqrt{n}%
:n\in\mathbb{N}\right\}  $, $2^{\mathbb{N}}=\left\{  2^{n}:n\in\mathbb{N}%
\right\}  $, and $h\mathbb{Z=}\left\{  hn:n\in\mathbb{Z}\right\}  $.
\end{example}

\section{Periodicity}

In the following we propose a new periodicity notion which does not oblige the
time scale to be closed under the operation $t\pm P$ for a fixed $P>0$ or to
be unbounded.

\begin{definition}
[Periodicity in shifts]\label{new per}Let $\mathbb{T}$ be a time scale with
the shift operators $\delta_{\pm}$ associated with the initial point $t_{0}%
\in\mathbb{T}^{\ast}$. The time scale $\mathbb{T}$ is said to be\emph{
\emph{periodic in} shifts }$\delta_{\pm}$ if there exists a $p\in(t_{0}%
,\infty)_{\mathbb{T}^{\ast}}$ such that $(p,t)\in\mathcal{D}_{\mp}$ for all
$t\in\mathbb{T}^{\ast}$. Furthermore, if%
\[
P:=\inf\left\{  p\in(t_{0},\infty)_{\mathbb{T}^{\ast}}:(p,t)\in\mathcal{D}%
_{\mp}\text{ for all }t\in\mathbb{T}^{\ast}\right\}  \neq t_{0},
\]
then $P$ is called the \emph{period} of the time scale $\mathbb{T}$.
\end{definition}

The following example indicates that a time scale, periodic in shifts, does
not have to satisfy (\ref{addition}). That is, a time scale periodic in shifts
may be bounded.

\begin{example}
\label{ex new per}The following time scales are not periodic in the sense of
Definition \ref{oldperts} but periodic with respect to the notion of shift
operators given in Definition \ref{new per}.

\begin{enumerate}
\item $\mathbb{T}_{1}\mathbb{=}\left\{  \pm n^{2}:n\in\mathbb{Z}\right\}  $,
$\delta_{\pm}(P,t)=\left\{
\begin{array}
[c]{ll}%
\left(  \sqrt{t}\pm\sqrt{P}\right)  ^{2} & \text{if }t>0\\
\pm P & \text{if }t=0\\
-\left(  \sqrt{-t}\pm\sqrt{P}\right)  ^{2} & \text{if }t<0
\end{array}
\right.  $, $P=1$, $t_{0}=0,$

\item $\mathbb{T}_{2}\mathbb{=}\overline{q^{\mathbb{Z}}}$, $\delta_{\pm
}(P,t)=P^{\pm1}t$, $P=q$, $t_{0}=1,$

\item $\mathbb{T}_{3}\mathbb{=}\overline{\mathbb{\cup}_{n\in\mathbb{Z}}\left[
2^{2n},2^{2n+1}\right]  }$, $\delta_{\pm}(P,t)=P^{\pm1}t$, $P=4$, $t_{0}=1,$

\item $\mathbb{T}_{4}\mathbb{=}\left\{  \frac{q^{n}}{1+q^{n}}:q>1\text{ is
constant and }n\in\mathbb{Z}\right\}  \cup\left\{  0,1\right\}  $,
\[
\delta_{\pm}(P,t)=\dfrac{q^{^{\left(  \frac{\ln\left(  \frac{t}{1-t}\right)
\pm\ln\left(  \frac{P}{1-P}\right)  }{\ln q}\right)  }}}{1+q^{\left(
\frac{\ln\left(  \frac{t}{1-t}\right)  \pm\ln\left(  \frac{P}{1-P}\right)
}{\ln q}\right)  }},\ \ P=\frac{q}{1+q}.
\]

\end{enumerate}
\end{example}

Notice that the time scale $\mathbb{T}_{4}$ in Example \ref{ex new per} is
bounded above and below and $\mathbb{T}_{4}^{\ast}=\left\{  \frac{q^{n}%
}{1+q^{n}}:q>1\text{ is constant and }n\in\mathbb{Z}\right\}  $.

\begin{remark}
Let $\mathbb{T}$ be a time scale that is periodic in shifts with the period
$P$. Thus, by P.4 of Definition \ref{shift} the mapping $\delta_{+}%
^{P}:\mathbb{T}^{\ast}\rightarrow\mathbb{T}^{\ast}$ defined by $\delta_{+}%
^{P}(t)=\delta_{+}(P,t)$ is surjective. On the other hand, we know by P.1 of
Definition \ref{shift} that shift operators $\delta_{\pm}$ are strictly
increasing in their second arguments. That is, the mapping $\delta_{+}%
^{P}(t):=\delta_{+}(P,t)$ is injective. Hence, $\delta_{+}^{P}$ is an
invertible mapping with the inverse $\left(  \delta_{+}^{P}\right)
^{-1}=\delta_{-}^{P}$ defined by $\delta_{-}^{P}(t):=\delta_{-}(P,t)$.
\end{remark}

In next two results, we suppose that $\mathbb{T}$ is a periodic time scale in
shifts $\delta_{\pm}$ with period $P$ and show that the operators $\delta
_{\pm}^{P}:\mathbb{T}^{\ast}\rightarrow\mathbb{T}^{\ast}$ are commutative with
the forward jump operator $\sigma:\mathbb{T\rightarrow T}$ given by
\[
\sigma(t):=\inf\left\{  s\in\mathbb{T}:s>t\right\}  .
\]
That is,%
\begin{equation}
\left(  \delta_{\pm}^{P}\circ\sigma\right)  (t)=\left(  \sigma\circ\delta
_{\pm}^{P}\right)  (t)\text{ for all }t\in\mathbb{T}^{\ast}.
\label{commutativity}%
\end{equation}

\begin{lemma}
The mapping $\delta_{+}^{T}:\mathbb{T}^{\ast}\rightarrow\mathbb{T}^{\ast}$
preserves the structure of the points in $\mathbb{T}^{\ast}$. That is,%
\[
\sigma(\widehat{t})=\widehat{t}\text{ implies }\sigma(\delta_{+}(P,\widehat
{t}))=\delta_{+}(P,\widehat{t}).
\]%
\[
\sigma(\widehat{t})>\widehat{t}\text{ implies }\sigma(\delta_{+}(P,\widehat
{t})>\delta_{+}(P,\widehat{t}).
\]

\end{lemma}

\begin{proof}
By definition we have $\sigma(t)\geq t$ for all $t\in\mathbb{T}^{\ast}$. Thus,
by P.1%
\[
\delta_{+}(P,\sigma(t))\geq\delta_{+}(P,t).
\]
Since $\sigma(\delta_{+}(P,t))$ is the smallest element satisfying%
\[
\sigma(\delta_{+}(P,t))\geq\delta_{+}(P,t),
\]
we get%
\begin{equation}
\delta_{+}(P,\sigma(t))\geq\sigma(\delta_{+}(P,t))\text{ for all }%
t\in\mathbb{T}^{\ast}\text{.} \label{1}%
\end{equation}
If $\sigma(\widehat{t})=\widehat{t}$, then (\ref{1}) implies%
\[
\delta_{+}(P,\widehat{t})=\delta_{+}(P,\sigma(\widehat{t}))\geq\sigma
(\delta_{+}(P,\widehat{t})).
\]
That is,%
\[
\delta_{+}(P,\widehat{t})=\sigma(\delta_{+}(P,\widehat{t}))\text{ provided
}\sigma(\widehat{t})=\widehat{t}\text{.}%
\]
If $\sigma(\widehat{t})>\widehat{t}$, then by definition of $\sigma$ we have%
\begin{equation}
(\widehat{t},\sigma(\widehat{t}))_{\mathbb{T}^{\ast}}=(\widehat{t}%
,\sigma(\widehat{t}))_{\mathbb{T}^{\ast}}=\varnothing\label{1.0}%
\end{equation}
and by P.1%
\[
\delta_{+}(P,\sigma(\widehat{t}))>\delta_{+}(P,\widehat{t}).
\]
Suppose contrary that $\delta_{+}(P,\widehat{t})$ is right dense, i.e.,
$\sigma(\delta_{+}(P,\widehat{t}))=\delta_{+}(P,\widehat{t})$. This along with
(\ref{1}) implies%
\[
(\delta_{+}(P,\widehat{t}),\delta_{+}(P,\sigma(\widehat{t})))_{\mathbb{T}%
^{\ast}}\neq\varnothing\text{.}%
\]
Pick one element $s\in(\delta_{+}(P,\widehat{t}),\delta_{+}(P,\sigma
(\widehat{t})))_{\mathbb{T}^{\ast}}$. Since $\delta_{+}(P,t)$ is strictly
increasing in $t$ and invertible there should be an element $t\in(\widehat
{t},\sigma(\widehat{t}))_{\mathbb{T}^{\ast}}$ such that $\delta_{+}(P,t)=s$.
This contradicts (\ref{1.0}). Hence, $\delta_{+}(P,\widehat{t})$ must be right
scattered, i.e., $\sigma(\delta_{+}(P,\widehat{t}))>\delta_{+}(P,\widehat{t}%
)$. The proof is complete.
\end{proof}

\begin{corollary}
\label{Cor 1} We have%
\begin{equation}
\delta_{+}(P,\sigma(t))=\sigma(\delta_{+}(P,t))\text{ for all }t\in
\mathbb{T}^{\ast}\text{.} \label{sigma delta1}%
\end{equation}
Thus,%
\begin{equation}
\delta_{-}(P,\sigma(t))=\sigma(\delta_{-}(P,t))\text{ for all }t\in
\mathbb{T}^{\ast}\text{.} \label{sigma delta2}%
\end{equation}

\end{corollary}

\begin{proof}
The equality (\ref{sigma delta1}) can be obtained as we did in the proof of
preceding lemma. By (\ref{sigma delta1}) we have%
\[
\delta_{+}(P,\sigma(s))=\sigma(\delta_{+}(P,s))\text{ for all }s\in
\mathbb{T}^{\ast}\text{.}%
\]
Substituting $s=\delta_{-}(P,t)$ we obtain%
\[
\delta_{+}(P,\sigma(\delta_{-}(P,t)))=\sigma(\delta_{+}(P,\delta
_{-}(P,t)))=\sigma(t)\text{.}%
\]
This and (iii) of Lemma \ref{lem pro} imply%
\[
\sigma(\delta_{-}(P,t))=\delta_{-}(P,\sigma(t))\text{ for all }t\in
\mathbb{T}^{\ast}.
\]
The proof is complete.
\end{proof}

Observe that (\ref{sigma delta1}) along with (\ref{sigma delta2}) yields
(\ref{commutativity}).

\begin{definition}
[Periodic function in shifts $\delta_{\pm}$]\label{def periodic in shift}Let
$\mathbb{T}$ be a time scale that is periodic in shifts $\delta_{\pm}$ with
the period $P$. We say that a real valued function $f$ defined on
$\mathbb{T}^{\ast}$ is\emph{ periodic in shifts }$\delta_{\pm}$ if there
exists a $T\in\lbrack P,\infty)_{\mathbb{T}^{\ast}}$ such that
\begin{equation}
\left(  T,t\right)  \in\mathcal{D}_{\pm}\text{ and }f(\delta_{\pm}%
^{T}(t))=f(t)\text{ for all }t\in\mathbb{T}^{\ast}, \label{t.1}%
\end{equation}
where $\delta_{\pm}^{T}(t):=\delta_{\pm}(T,t)$. The smallest number
$T\in\lbrack P,\infty)_{\mathbb{T}^{\ast}}$ such that (\ref{t.1}) holds is
called the period of $f$.
\end{definition}

\begin{example}
By Definition \ref{new per} we know that the set of reals $\mathbb{R}$ is
periodic in shifts $\delta_{\pm}$ defined by (\ref{rs1}-\ref{rs2}) associated
with the initial point $t_{0}=1$. The function%
\[
f(t)=\sin\left(  \frac{\ln\left\vert t\right\vert }{\ln\left(  1/2\right)
}\pi\right)  \text{, }t\in\mathbb{R}^{\ast}:=\mathbb{R}-\left\{  0\right\}
\]
is periodic in shifts $\delta_{\pm}$ defined by (\ref{rs1}-\ref{rs2}) with the
period $T=4$ since%
\begin{align*}
f\left(  \delta_{\pm}(T,t)\right)   &  =\left\{
\begin{array}
[c]{cc}%
f\left(  t4^{\pm1}\right)  & \text{if }t\geq0\\
f\left(  t/4^{\pm1}\right)  & \text{if }t<0\text{ }%
\end{array}
\right. \\
&  =\sin\left(  \frac{\ln\left\vert t\right\vert \pm2\ln\left(  1/2\right)
}{\ln\left(  1/2\right)  }\pi\right) \\
&  =\sin\left(  \frac{\ln\left\vert t\right\vert }{\ln\left(  1/2\right)  }%
\pi\pm2\pi\right) \\
&  =\sin\left(  \frac{\ln\left\vert t\right\vert }{\ln\left(  1/2\right)  }%
\pi\right) \\
&  =f(t)
\end{align*}
for all $t\in\mathbb{R}^{\ast}$ (see Figure 1).

\begin{figure}[h]
%Requires \usepackage{graphicx}
\includegraphics[width=3.5in]{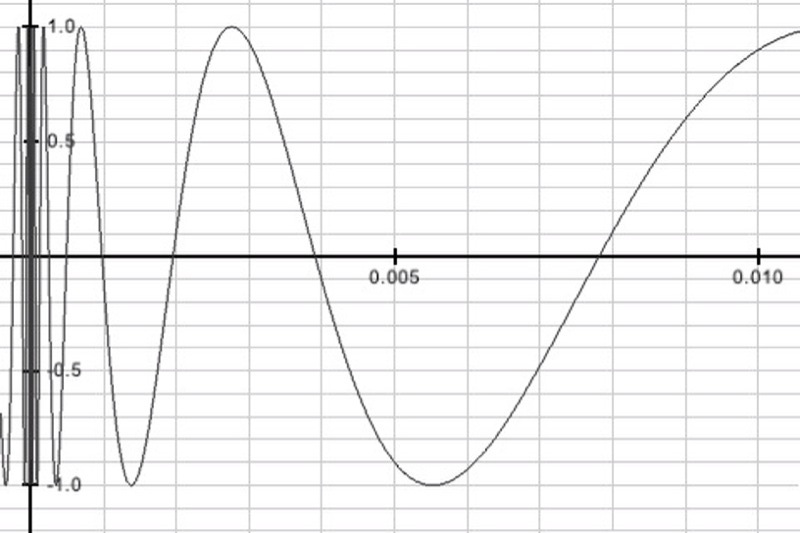}\newline\caption{Graph of $f(t)=\sin
\left(  \frac{\ln\left\vert t\right\vert }{\ln\left(  1/2\right)  }\pi\right)
$}%
\end{figure}
\end{example}

\begin{example}
The time scale $\overline{q^{\mathbb{Z}}}=\left\{  q^{n}:n\in\mathbb{Z}\text{
and }q>1\right\}  \cup\left\{  0\right\}  $ is periodic in shifts $\delta
_{\pm}(P,t)=P^{\pm1}t$ with the period $P=q$. The function $f$ defined by%
\begin{equation}
f(t)=\left(  -1\right)  ^{\dfrac{\ln t}{\ln q}},\ \ t\in q^{\mathbb{Z}}
\label{f}%
\end{equation}
is periodic in shifts $\delta_{\pm}$ with the period $T=q^{2}$ since
$\delta_{+}(q^{2},t)\in\overline{q^{\mathbb{Z}}}^{\ast}=q^{\mathbb{Z}}$ and%
\[
f\left(  \delta_{\pm}(q^{2},t)\right)  =\left(  -1\right)  ^{\dfrac{\ln t}{\ln
q}\pm2}=\left(  -1\right)  ^{\dfrac{\ln t}{\ln q}}=f(t)
\]
for all $t\in q^{\mathbb{Z}}$. However, $f$ is not periodic in the sense of
Definition \ref{old per} since there is no any positive number $T$ so that
$f(t\pm T)=f(t)$ holds.
\end{example}

In the following, we introduce $\Delta$-periodic function in shifts. For a
detailed information on $\Delta-$derivative and $\Delta$-integration we refer
to \cite{book} and \cite{book2}.

\begin{definition}
[$\Delta-$periodic function in shifts $\delta_{\pm}$]%
\label{def delta periodic in shifts}Let $\mathbb{T}$ be a time scale that is
periodic in shifts $\delta_{\pm}$ with period $P$. We say that a real valued
function $f$ defined on $\mathbb{T}^{\ast}$ is\emph{ }$\Delta-$\emph{periodic
in shifts }$\delta_{\pm}$ if there exists a $T\in\lbrack P,\infty
)_{\mathbb{T}^{\ast}}$ such that
\begin{equation}
\left(  T,t\right)  \in\mathcal{D}_{\pm}\text{ for all }t\in\mathbb{T}^{\ast},
\label{t.3}%
\end{equation}%
\begin{equation}
\text{the shifts }\delta_{\pm}^{T}\text{ are }\Delta-\text{differentiable with
rd-continuous derivatives,} \label{t3.1}%
\end{equation}
and%
\begin{equation}
f(\delta_{\pm}^{T}\left(  t\right)  )\delta_{\pm}^{\Delta T}(t)=f(t)
\label{t4}%
\end{equation}
for all $t\in\mathbb{T}^{\ast}$, where $\delta_{\pm}^{T}(t):=\delta_{\pm
}(T,t)$. The smallest number $T\in\lbrack P,\infty)_{\mathbb{T}^{\ast}}$ such
that (\ref{t.3}-\ref{t4}) hold is called the period of $f$.
\end{definition}

Notice that Definition \ref{def periodic in shift} and Definition
\ref{def delta periodic in shifts} give the classic periodicity definition
(i.e. Definition \ref{old per}) on time scales whenever $\delta_{\pm}%
^{T}(t)=t\pm T$ are the shifts satisfying the assumptions of Definition
\ref{def periodic in shift} and Definition \ref{def delta periodic in shifts}.

\begin{example}
The real valued function $g(t)=1/t$ defined on $2^{\mathbb{Z}}=\{2^{n}%
:n\in\mathbb{Z}\}$ is $\Delta-$\emph{periodic in shifts} $\delta_{\pm
}(T,t)=T^{\pm1}t$ with the period $T=2$ since%
\[
f\left(  \delta_{\pm}(2,t)\right)  \delta_{\pm}^{\Delta}(2,t)=\frac{1}%
{2^{\pm1}t}2^{\pm1}=\frac{1}{t}=f(t).
\]
The following result is essential for the proof of next theorem
\end{example}

\begin{theorem}
[Substitution]\label{thm2.2} \cite[Theorem 1.98]{book} Assume $\nu
:\mathbb{T}\rightarrow\mathbb{R}$ is strictly increasing and $\tilde
{\mathbb{T}}:=\nu(\mathbb{T})$ is a time scale. If $f:\mathbb{T}%
\rightarrow\mathbb{R}$ is an rd-continuous function and $\nu$ is
differentiable with rd-continuous derivative, then for $a,b\in\mathbb{T}$,
\begin{equation}
\int_{a}^{b}\!g(s)\nu^{\Delta}(s)\,\Delta s=\int_{\nu(a)}^{\nu(b)}g(\nu
^{-1}(s))\,\tilde{\Delta}s. \label{substitute}%
\end{equation}

\end{theorem}

\begin{theorem}
Let $\mathbb{T}$ be a time scale that is periodic in shifts $\delta_{\pm}$
with period $P\in(t_{0},\infty)_{\mathbb{T}^{\ast}}$ and $f$ a\emph{ }%
$\Delta-$periodic function in shifts\emph{ }$\delta_{\pm}$ with the period
$T\in\lbrack P,\infty)_{\mathbb{T}^{\ast}}$. Suppose that $f\in C_{rd}%
(\mathbb{T})$, then%
\[
\int_{t_{0}}^{t}f(s)\Delta s=\int_{\delta_{\pm}^{T}(t_{0})}^{\delta_{\pm}%
^{T}(t)}f(s)\Delta s.
\]

\end{theorem}

\begin{proof}
Substituting $v(s)=$ $\delta_{+}^{T}(s)$ and $g(s)=f(\delta_{+}^{T}\left(
s\right)  )$ in (\ref{substitute}) and taking (\ref{t4}) into account we have%
\begin{align*}
\int_{\delta_{+}^{T}(t_{0})}^{\delta_{+}^{T}(t)}f(s)\Delta s  &  =\int
_{\nu(t_{0})}^{\nu(t)}g(\nu^{-1}(s))\Delta s\\
&  =\int_{t_{0}}^{t}g(s)\nu^{\Delta}(s)\Delta s\\
&  =\int_{t_{0}}^{t}f(\delta_{+}^{T}\left(  s\right)  )\delta_{+}^{\Delta
T}(t)\Delta s\\
&  =\int_{t_{0}}^{t}f(s)\Delta s.
\end{align*}
The equality%
\[
\int_{\delta_{-}^{T}(t_{0})}^{\delta_{-}^{T}(t)}f(s)\Delta s=\int_{t_{0}}%
^{t}f(s)\Delta s
\]
can be obtained similarly. The proof is complete.
\end{proof}

\begin{center}
{\large Acknowledgement}
\end{center}

I would like to thank reviewers for their valuable comments improving this study.


\begin{thebibliography}{xxx}                                                                                               %


\bibitem {adivar}
ADIVAR, M:
\textit{Function bounds for solutions of Volterra integro dynamic equations on time scales}, 
E. J. Qualitative Theory of Diff. Equ., \textbf{7}, (2010), 1--22.

\bibitem {bams}
ADIVAR, M.---RAFFOUL, Y. N.: 
\textit{Existence of resolvent for Volterra integral equations on time scales}, 
Bull. of Aust. Math. Soc., \textbf{82 (1)}, (2010), 139--155.

\bibitem {adivar1}ADIVAR, M.---RAFFOUL, Y. N.: 
\textit{Stability and periodicity in dynamic delay equations}, 
Computers and Mathematics with Applications, \textbf{58 (2)}, (2009), 264--272.

\bibitem {adivar2}ADIVAR, M.---RAFFOUL, Y. N.:
\textit{Existence results for periodic solutions of integro-dynamic equations on time scales}, 
Annali di Matematica ed Pure Applicata, \textbf{188 (4)}, (2009), 543--559.

\bibitem {atici}
ATICI, F. M.---GUSEINOV, G. SH.---KAYMAK\c{C}ALAN, B.:
\textit{Stability criteria for dynamic equations on time scales with periodic coefficients}.
In: Proceedings of the International Conference on Dynamic Systems and Applications III, (3) (1999), 43--48.

\bibitem {Bi&Bohner}BI, L.---BOHNER, M.---FAN, M.:
\textit{Periodic solutions of functional dynamic equations with infinite delay}, 
Nonlinear Anal., \textbf{68(5)}, (2008), 1226--1245.

\bibitem {book}BOHNER, M.---PETERSON, A. C.:
\textit{Dynamic equations on time scales, An introduction with applications},
Birkh\"{a}user Boston Inc., Boston, MA, 2001.

\bibitem {book2}BOHNER, M.---PETERSON, A. C.: 
\textit{Advances in Dynamic equations on time scales}, 
Birkh\"{a}user Boston Inc., Boston, MA, 2003.

\bibitem {Kaufmann 2}KAUFMANN E.---RAFFOUL, Y. N.: 
\textit{Periodic solutions for a neutral nonlinear dynamical equations on time scale}, 
J. Math. Anal. Appl., \textbf{319, (1)}, (2006), 315--325.

\bibitem {Kaufmann&raffoul}KAUFMANN E.---RAFFOUL, Y. N.:, 
\textit{Periodicity and stability in neutral nonlinear dynamic equations with functional delay on a
time scale}, 
Electron. J. Differential Equations, \textbf{No. 27}, (2007), 1--12.

\bibitem {li}LIUA, XI-LAN---LIB, WAN-TONG: 
\textit{Periodic solutions for dynamic equations on time scales}, 
Nonlinear Analysis, \textbf{67}, (2007), 1457--1463.
\end{thebibliography}
\end{document}